\documentclass{ifacconf}
\usepackage{amsfonts,amssymb,mathrsfs,amsmath}
\usepackage{graphicx} 
\usepackage{xcolor}
\usepackage{array}
\usepackage{natbib}   
\newtheorem{Rem}{Remark}
\newtheorem{Def}{Definition}
\newcommand{\N}{\mathbb{N}}
\newcommand{\R}{\mathbb{R}}
\renewcommand{\S}{\mathbf{S}}
\begin{document}
\begin{frontmatter}

\title{Symplectic Methods in Deep Learning\thanksref{footnoteinfo1}\thanksref{footnoteinfo2}} 

\thanks[footnoteinfo1]{This work has been supported by Deutsche Forschungsgemeinschaft (DFG), Grant No. OB 368/5-1, AOBJ: 692093}

\thanks[footnoteinfo2]{\copyright 2024 the
authors. This work has been accepted to MTNS 2024 for publication under a Creative Commons Licence CC-BY-NC-ND}

\author[First]{Sofya Maslovskaya} 
\author[Second]{Sina Ober-Blöbaum} 

\address[First]{Institute of Mathematics, Paderborn University (e-mail: sofyam@math.uni-paderborn.de).}
\address[Second]{Institute of Mathematics, Paderborn University (e-mail: sinaober@math.uni-paderborn.de).}

\begin{abstract}                
Deep learning is widely used in tasks including image recognition and generation, in learning dynamical systems from data and many more. It is important to construct learning architectures with theoretical guarantees to permit safety in the applications. There has been considerable progress in this direction lately. In particular, symplectic networks were shown to have the non vanishing gradient property, essential for numerical stability. On the other hand, architectures based on higher order numerical methods were shown to be efficient in many tasks where the learned function has an underlying dynamical structure. In this work we construct symplectic networks based on higher order explicit methods with non vanishing gradient property and test their efficiency on various examples.
\end{abstract}

\begin{keyword}
Deep neural network, non-autonomous Hamiltonian system, symplectic integrator, partitioned Runge-Kutta method, higher order method, learning Hamiltonian systems.
\end{keyword}

\end{frontmatter}

\section{Introduction}
Deep learning proved to be efficient in different learning tasks, especially in image classification and learning of dynamical systems. In this context, design of networks with theoretical guarantees is of special interest for safety in applications. There are properties, which are crucial for the learning process. Vanishing and exploding gradients are the problems of this kind and known to cause instabilities in the optimization process and failure in learning (\cite{Bengio1994}). 

It is shown in different studies that considering neural networks based on dynamical systems can be helpful to design networks and gradient computation algorithms enjoying theoretically proven properties, such as stability or structure preservation (\cite{Chen2018,yan2019,CELLEDONI2021}). In particular, it has been recently discovered in (\cite{Galimberti2023}) that Hamiltonian networks possess the non vanishing gradient property. The network architecture is given by a symplectic Euler discretization of a certain separable Hamiltonian system, which plays the role of the activation function. Even though, Hamiltonian networks are of a very particular structure, the universal approximation properties is shown in (\cite{Zakwan2023}) based on the feature augmentation (\cite{Dupont2019}). 

Recent results (\cite{Giesecke2021,Benning2021}) show advantages in using higher order discretizations in network architecture. In particular, one can reach much better approximation properties in learning dynamical systems when the network represents a higher order numerical method according to (\cite{Matsubara2023}). This is why we extend the existing Hamiltonian network from (\cite{Galimberti2023}) to a new class of networks based on higher order symplectic partitioned Runge-Kutta (SPRK) discretizations. 

%
%
%
%
 The contribution of this work is as follows. We introduce a new class of neural networks based on explicit SPRK methods applied to a non autonomous separable Hamiltonian system. For this, we obtain the condition on a PRK method to be symplectic, explicit and of higher order. The new construction is built upon the classical conditions for autonomous Hamiltonian systems. The non vanishing gradient property of the constructed network follows from the symplecticity of the used numerical methods. In addition, we show the universal approximation property of the network. Furthermore, the networks are based on explicit numerical methods, which implies low computational cost with respect to implicit analogues. Learning problems formulated with the network architecture based on a higher order method have the property to be a higher order approximation of a continuous optimal control problem. This property can be advantageous for the further analysis and permits to better understand the nature of the transformation which is happening in the neural network. We test the obtained network class on different learning tasks and compare it with the Hamiltonian network from (\cite{Galimberti2023}). The learning examples include the image classification and the learning of an autonomous Hamiltonian system.  
\section{Preliminaries}

\subsection{Deep learning}
The main  goal of machine learning is to provide an approximation for an unknown function $F : X \rightarrow Y$ by generating a function $\hat F$ in such a way that $\hat F$ approximates the given training data $\{x_i,  y_i = F(x_i) \}^{n}_{i = 1}$, where the arguments $x_i$ are called features and the function values $y_i$ are called labels. In addition, it must generalize to unknown data $\{\tilde x_i,  \tilde y_i = F(\tilde x_i) \}^{m}_{i = 1}$. Deep learning is a variant of machine learning where the approximation is realized by a neural network with a large number of layers. A neural network represents a composition of parameterized nonlinear functions $\sigma(W_j (x_{i})_j + \beta_j)$. The function $\sigma$ is called activation function. The propagation of feature $x_i$ through the network can be represented in the following way.
\begin{multline*}
x_i \rightarrow \sigma(W_0 x_i + \beta_0) \rightarrow \sigma ( W_1\sigma(W_0 x_i + \beta_0) + \beta_1) \rightarrow \cdots \\ \cdots \rightarrow \sigma( W_{N-1}\sigma(\cdots) + \beta_{N-1}) ) \rightarrow\hat F_{W,\beta}( x_i).
\end{multline*}
In order to find the parameters $\{W_k,\beta_k\}_{0 \leq k \leq N-1}$ for which $\hat F_{W,\beta}$ approximates $F$ in the best possible way, an optimization problem is formulated. The objective function 
measures how far $\hat F_{W,\beta}(x_i)$ is from $y_i$. All the features of the training data are composed into a vector $z = (x_1, \dots, x_n)$ and labels into $y = (y_1, \dots, y_m)$. The optimization problem is defined as follows
\begin{equation} \label{eq:DL:optimal.problem}
\begin{aligned}
& \min_{\{ W, \beta \}} \mathcal{J}(W,\beta) = L(\Phi(z_N), y) + R(W, \beta) \\
&  z_{k+1}=\sigma (W_k z_{k} + \beta_k), \quad k = 0, \dots, N-1, \\
& z_0 = z, 
\end{aligned}
\end{equation}
where $\Phi(\cdot)$ is the output activation function defined on the space of features with its image in the space of labels, $L(\cdot, \cdot)$ is a loss function which measures the distance between the output of the network 
and the labels and $R$ is the regularization term.
			
\subsubsection{Gradients} 			
Typically, variants of stochastic gradient descent are employed to solve \eqref{eq:DL:optimal.problem}. Gradient based methods find the optimal solution using the gradient descent until a critical point with zero-gradient is found
\begin{equation} \label{eq:zero.gradient}
    \partial_{u} \left( L(u) + R(u)  \right) = 0, 
\end{equation}

with $u$ the collection of all optimization parameters $\{W_k,\beta_k\}_{0 \leq k \leq N-1}$. The next parameter value is found by
$$ u_{j+1} = u_j  + \gamma \, \partial_{u} \left( L(u_j) + R(u_j)  \right),$$
$\gamma$ is a learning rate. The calculation of gradients is done using the backpropagation algorithm. The relevance of backward propagation is a direct consequence of the composition structure of the network. The gradients are computed by the following formula
$$\frac{\partial L}{\partial u_k} = \frac{\partial L}{\partial z_N} \left( \prod_{i = k+1}^{N-1} \frac{\partial z_{i+1}}{\partial z_{i}}  \right) \frac{\partial z_{k+1}}{\partial u_k}. $$
The main idea of the backpropagation is the application of reverse-mode auto-differentiation. Theoretically, deeper networks are expected to provide a better approximation, but in practice we usually observe the phenomenon of vanishing or exploding gradient which are due to the error accumulation in the computation of gradients. The vanishing gradient issue occurs in case of the smallness of $\| \prod_{i = k+1}^{N-1} \frac{\partial z_{i+1}}{\partial z_{i}}  \|$.

\subsubsection{Universal approximation property}
One of the central requirements for a network is its universal approximation property. This property characterizes the ability of the given network to approximate any function in a given function class as accurately as possible. More precisely, we use the following definition as in (\cite{Zakwan2023})
\begin{Def}[UAP] \label{def:UAP}
    A neural network $\hat F_{W,\beta}$ parameterized by $\{W_k,\beta_k\}_{0 \leq k \leq N-1}$ has a universal approximation property (UAP) on the set of continuous functions defined on a compact domain  $\Omega \in \R^{n_1}$ with values in $\R^{n_2}$ if for any such function $f$ 
    and any $\varepsilon > 0$, there exist a number of layers $N$ and parameters $\{W_k,\beta_k\}_{0 \leq k \leq N-1}$  such that
    $$ \sup_{x \in\Omega} \| f(x) - \hat F_{W,\beta}(x) \| \leq \varepsilon.$$
\end{Def}
\subsubsection{Neural ODEs}
Neural ordinary differential equations (neural ODEs) come into play when we assume a particular structure of the activation function $\sigma$, to be a discretization of an ODE. In this case
\begin{equation} \label{eq:neural.ODE}
z_{k+1}=\sigma (W_k z_{k} + \beta_k)    
\end{equation}
approximates
$$\dot z(t) = f(W(t)z(t) + \beta(t)),$$
for some vector field $f$, where $z(t), W(t), \beta(t)$ are absolutely continuous functions of time of an appropriate dimension. Then, $z_k, W_k, \beta_k$ in \eqref{eq:neural.ODE} approximate $z(t_k)$, $W(t_k), \beta(t_k)$ at the discrete grid of time points $t_1, \dots, t_N$. This approach to activation function design permits to construct neural networks preserving the desired properties of the flow of a dynamical system.  
In particular, in this work we consider neural networks based on symplectic discretizations of Hamiltonian systems to guarantee the non vanishing gradient property.

\subsection{Symplectic integration of Hamiltonian systems}

Let us consider a Hamiltonian system defined on a cotangent space $T^*M$ of a manifold $M$ by a time dependent Hamiltonian $H$. In canonical coordinates $(q,p) \in T^*M$, the system is given by
\begin{equation} \label{eq:H.syst}
\dot{q} = \frac{\partial H}{\partial p}(q,p,t), \quad \dot{p} = -\frac{\partial H}{\partial q}(q,p,t). 
\end{equation}
The flow $(q(t),p(t))$ preserves the canonical symplectic form $dq \wedge dp$. In fixed coordinates, it can be also written in a matrix form
$$J = \begin{pmatrix}
    0 & I \\
    -I & 0
\end{pmatrix}.$$
Let $\phi_t$ be the flow of \eqref{eq:H.syst} from a given initial point $z_0 = (q_0, p_0)$. The preservation of the symplectic form can be expressed as follows
$$
\left(\frac{\partial \phi_t}{\partial z_0}\right)^\top J \left(\frac{\partial \phi_t}{\partial z_0}\right) = J.
$$
This implies, in particular, that the Jacobian of the flow with respect to the initial data satisfies
$$
\left\|\frac{\partial \phi_t}{\partial z_0} \right\| \geq 1.
$$
We would like to have this property on the level of the neural network. This can be done by using symplectic methods to approximate the flow of the Hamiltonian system. More concretely, 
let us consider a one step numerical method $(q_n, p_n) \rightarrow (q_{n+1}, p_{n+1})$ applied to \eqref{eq:H.syst}
			\begin{equation*} 
				\begin{pmatrix}
					q_{n+1} \\ p_{n+1}
				\end{pmatrix} = \begin{pmatrix}
					q_n + h \varphi_q(q_n, p_n, q_{n+1}, p_{n+1}, t_n, t_{n+1}) \\
					p_n + h \varphi_p(q_n, p_n, q_{n+1}, p_{n+1}, t_n, t_{n+1})
				\end{pmatrix}.
			\end{equation*} 
The methods which preserve the symplectic form along the numerical flow are called symplectic (\cite{Hairer2013}). The preservation of the symplectic form in numerical flow means: $dq_{n+1} \wedge dp_{n+1} = dq_{n} \wedge dp_{n}$ or in coordinates
			\begin{equation} \label{eq:sympl.num.preservation}
				\begin{pmatrix}
					\frac{\partial q_{n+1}}{\partial q_{n}} & \frac{\partial q_{n+1}}{\partial p_{n}} \\
					\frac{\partial p_{n+1}}{\partial q_{n}} & \frac{\partial p_{n+1}}{\partial p_{n}}
				\end{pmatrix}^\top J \begin{pmatrix}
					\frac{\partial q_{n+1}}{\partial q_{n}} & \frac{\partial q_{n+1}}{\partial p_{n}} \\
					\frac{\partial p_{n+1}}{\partial q_{n}} & \frac{\partial p_{n+1}}{\partial p_{n}}
				\end{pmatrix} = J.
    \end{equation}

\subsection{Hamiltonian neural network}
It is suggested in (\cite{Galimberti2023}) to use the symplectic Euler method in combination with a separable Hamiltonian  system for the construction of a neural network with non vanishing gradient property. The key property to ensure the non vanishing gradient in this case is given by \eqref{eq:sympl.num.preservation}. Following (\cite{Galimberti2023}), let us consider the following separable Hamiltonian
\begin{equation} \label{eq:network:Hamiltonian}
    H(z,t) = \left( ( 1, \cdots, 1) \, \tilde \sigma(W z + \beta)	 + \eta^\top z \right) 
\end{equation} 

		with parameters
		$$ W(t) = \begin{pmatrix}
		W_1(t) & 0 \\
	     0 & W_2(t)
		\end{pmatrix}, \, \beta(t) = \begin{pmatrix}
		\beta_1(t) \\
		\beta_2(t)
	\end{pmatrix}, \, \eta(t) = \begin{pmatrix}
	\eta_1(t) \\
	\eta_2(t)
\end{pmatrix}, $$
and $\tilde \sigma = (\tilde \sigma_1, \dots, \tilde \sigma_n)$ with $\tilde \sigma_i : \R \rightarrow \R$ for $i = 1, \dots, n$, which is related to the activation function $\sigma = (\sigma_1, \dots, \sigma_n)$ by $\tilde \sigma_i(\cdot)' = \sigma_i(\cdot)$. The corresponding Hamiltonian system for $z = (q,p)$ is defined as follows
\begin{equation}
\begin{aligned} \label{eq:Ham.syst}
	&\dot{q} = W_2^\top \sigma(W_2 p + \beta_2) + \eta_2,\\
	& \dot{p} = -W_1^\top \sigma(W_1 q + \beta_1) + \eta_1.
\end{aligned}
\end{equation}
Applying symplectic Euler to it leads to the following explicit scheme
\begin{equation} \label{eq:H.net}
\begin{aligned}
	& p_{j+1} = p_j - h \left( W_{1,j}^\top \sigma (W_{1,j} q_j + \beta_{1,j}) + \eta_{1,j}\right), \\
	& q_{j+1} = q_j +h \left( W_{2,j}^\top \sigma (W_{2,j} p_{j+1} + \beta_{2,j}) + \eta_{2,j}\right).
\end{aligned}    
\end{equation}

It is shown in (\cite{Galimberti2023}), that the neural network of the form \eqref{eq:H.net} called Hamiltonian Net is symplectic. 
Preservation of the symplectic form $J$ by a forward pass from $z_n$ to $z_{n+1}$ of a network implies the non-vanishing gradient property. This follows from 
(\cite{Galimberti2023})
\begin{equation} \label{eq:sympl.nonvanishing}
 \frac{\partial z_{n+1}}{\partial z_{n}}^\top J \frac{\partial z_{n+1}}{\partial z_{n}} = J \ \Rightarrow \ \left\|\frac{\partial z_{n+1}}{\partial z_{n}} \right\|  \geq 1, 
 \end{equation}
 which implies
 \begin{equation}
   \left\| \prod_{i = k+1}^{N-1} \frac{\partial z_{i+1}}{\partial z_{i}} \right \| \geq 1. 
\end{equation} 
%
Moreover, it is shown in (\cite{Zakwan2023}) that using the feature space augmentation, it is possible to have the universal approximation theorem for Hamiltonian Net. The feature augmentation is done by doubling the dimension of the feature space. If the features in the training data are given by $z_0 \in \R^n$, then the new augmented feature space is $\R^{2n}$. Hamiltonian Net \eqref{eq:H.net} is then considered in $\R^{2n}$ with the new variables $(q_n, p_n)$ and the initial data $(q_0, p_0) = (0, z_0) \in \R^{2n}$. 
In this case the following theorem holds.
\begin{thm}[\cite{Zakwan2023}] \label{th:h.Net.UAP1}
Given a compact $\Omega \in \R^{n}$ and $\xi \in \Omega$, the Hamiltonian network defined by \eqref{eq:H.net} with $(q_0, p_0) = (0, \xi)$, $q_n, p_n \in \R^{n}$ for $n = 0, \dots, N$ and the output $q_N$, that is projection of  $(q_N, p_N)$ to the first component, satisfies the universal approximation property (Definition~\ref{def:UAP}) for functions $f : \Omega \rightarrow \R^{n}$.
\end{thm}
An adapted output layer permits to generalize the UAP to general functions between spaces of different dimensions.
\begin{thm}[\cite{Zakwan2023}] \label{th:h.Net.UAP2}
If a network satisfies the UAP for all continuous functions $f : \Omega \subset \R^{n} \rightarrow \R^{n}$, then there exist parameters $ W_0,b_0$  of appropriate dimensions, such that the network composed with an output activation function $h(q,p) = W_0 (q,p)^\top + b_0$ satisfies the universal approximation property for any continuous function on $\Omega$.
\end{thm}

\section{Higher order Hamiltonian network}
Let us now construct a new class of higher order Hamiltonian networks based on higher order symplectic methods. The constructed new class of networks will also have the advantages of the Hamiltonian network from the previous section.
\subsection{Symplectic partitioned Runge-Kutta methods}
The use of implicit methods is not desirable in neural network architecture because it requires a high number of computations and might course additional problems in the gradient calculations. This is why our goal is to use explicit symplectic methods. The most appropriate class of numerical methods for separable systems is
the Partitioned Runge-Kutta (PRK) method. It can be defined using Butcher tables as follows (\cite{SanzSerna2016})
			 $$
				\begin{array}
					{c|cccc}
					c_1 & a_{11} & \dots & a_{1s} \\
					\vdots & \vdots \\
					c_s & a_{s1} & \dots & a_{ss} \\
					\hline
					& b_1 & \dots & b_s 
				\end{array} \qquad \begin{array}
				{c|cccc}
				C_1 & A_{11} & \dots & A_{1s} \\
				\vdots & \vdots \\
				C_s & A_{s1} & \dots & A_{ss} \\
				\hline
				& B_1 & \dots & B_s 
			\end{array}
			$$
The next step $(q_{n+1}, p_{n+1})$ is obtained from $(q_n, p_n)$ by 
			\begin{equation} \label{general.PRK}
			\begin{aligned}
			& q_{n+1} = q_n+h \sum_{i=1}^{s} b_i k_{n,i}, \\
			& k_{n,i} = \frac{\partial H}{\partial p}(Q_{n,i},P_{n,i}, t_n + c_i h), \\
			& Q_{n,i} = q_n + h \sum_{i=1}^{s} a_{ij} k_{n,j}, \\
			& p_{n+1} = p_n+h \sum_{i=1}^{s} B_i \ell_{n,i}, \\
			& \ell_{n,i} = -\frac{\partial H}{\partial q}(Q_{n,i},P_{n,i}, t_n + C_i h), \\
			& P_{n,i} = p_n + h \sum_{i=1}^{s} A_{ij} \ell_{n,j}.
		\end{aligned}
            \end{equation}
  The following theorem gives the conditions for a PRK method applied to a separable Hamiltonian to be symplectic. The proof  is based on the proof presented in (\cite{abia1993}) for the autonomous Hamiltonian and the statement for the non autonomous case is mentioned in (\cite{Kalogiratou2011}) without proof.
  \begin{prop}
If $H(q,p,t)$ is separable in $(q,p)$,  i.e., $H(q,p,t) = H_1(q,t) + H_2(p,t)$, then the following condition on $A_{ij}, a_{i,j}, B_i, b_i$ implies symplecticity of a PRK method
		$$ b_i A_{ij} + B_j a_{i,j} - b_i B_j = 0, \qquad i,j = 1 \dots s.$$
\end{prop}
\begin{pf} $dq_{n+1} \wedge dp_{n+1} - dq_{n} \wedge dp_{n} = $ 
	\begin{multline*}
	 h \sum_{i=1}^{s} \underbrace{\left( b_i dQ_{n,i} \wedge dk_{n,i}  + B_i d\ell_{n,i} \wedge dP_{n,i}\right)}_{ \left( b_i \frac{\partial^2 H}{\partial q \partial p}(t+c_ih) - B_i \frac{\partial^2 H}{\partial q \partial p}(t+C_ih) \right)dQ_{n,i} \wedge dP_{n,i} } \\- h^2 \sum_{i,j = 1}^{s} \left( b_i A_{ij} + B_j a_{i,j} - b_i B_j \right)d\ell_{n,j} \wedge d k_{n,i}.
	\end{multline*} 
For a separable Hamiltonian, we have $\frac{\partial^2 H}{\partial q \partial p} = 0$. 
\end{pf}
\begin{Rem}
    Notice that the condition for separable Hamiltonian is different from the general case with the additional conditions $b_i = B_i$ and $c_i = C_i$. This permits the construction of explicit methods.
\end{Rem}
The condition above coincides with the condition for time invariant Hamiltonian systems. This permits to use classical form of explicit PRK methods for autonomous systems, such as methods constructed in (\cite{Forest1990}), defined by the following Butcher tables. 
			 \begin{equation} \label{eq:SPRK.butcher}
			    \begin{array}
				{c|cccc}
				c_1 & b_{1} & 0 & \dots & 0\\
				c_2 & b_{1} & b_{2} & \dots & 0\\
				\vdots & \vdots \\
				c_s & b_{1}& b_{2} & \dots & b_{s} \\
				\hline
				& b_1 & b_{2} & \dots & b_s 
			\end{array} \qquad \begin{array}
			{c|cccc}
			C_1 & 0 & 0 & \dots & 0\\
			C_2 & B_{1} & 0 & \dots & 0\\
			\vdots & \vdots \\
			C_s & B_{1}& B_{2} & \dots & 0 \\
			\hline
			& B_1 & B_{2} & \dots & B_s 
		\end{array} 
			 \end{equation}
The table above is completely determined by $b = (b_1, \dots, b_s), \ B = (B_1, \dots, B_s), \ c =  (c_1, \dots, c_s)$ and $C=(C_1, \dots, C_s)$. The corresponding explicit scheme is 
\begin{equation} \label{eq:expl.SPRK}
 \begin{aligned}
			& Q_{n,1} = q_n \\
			& P_{n,1} = p_n \\
                & \text{for }i=1 \dots s-1 \\
			& \ \ \ \ Q_{n,i+1} = Q_{n,i} + h b_{i} k_{n,i}(P_{n,i}, u(t_n + c_i h)) \\
			& \ \ \ \ P_{n,i+1} = P_{n,i} + h B_{i} \ell_{n,i+1}(Q_{n,i+1}, u(t_n + C_i h)) \\
			&  q_{n+1} = Q_{n,s}\\
			&  p_{n+1} = P_{n,s}\\
		\end{aligned}   
\end{equation}

Notice that \eqref{eq:SPRK.butcher} does not allow to consider the same time stages, i.e.\ $c = C$, as this would lead to a non consistent numerical scheme. 
 The restriction $c=C$ is typically used for the approximation of non autonomous Hamiltonian systems, where only potential energy is time-dependent (\cite{Forest1990,Hairer2013}). However, it is crucial  for the network expressivity to consider a Hamiltonian $H(q,p,t) = H_1(q,t) + H_2(p,t)$, where both $H_1$ and $H_1$ depend on time. This is required in the proof of Theorem~\ref{th:univ.approx}. 
%
Nevertheless, it is possible to use the order conditions on $b, B$ for autonomous systems in combination with a specific condition on $c, C$ to obtain higher order SPRK of the form \eqref{eq:SPRK.butcher}.
\begin{thm} \label{thm:order.nonaut}
    Assume that $b, a = (a_{ij}), B, A = (A_{ij})$ define PRK~\eqref{general.PRK} and satisfy the $p$-order condition for autonomous separable Hamiltonian system, i.e., for all $H(q,p) = H_1(q)+H_2(q)$, with $p\leq 4$ and in addition $c_i = \sum_{j = 1}^s A_{ij}$ and $C_i = \sum_{j = 1}^s a_{ij}$. Then the PRK method defined by $b, a, c, B, A, C$  has order $p$ for all non autonomous separable Hamiltonian systems.
\end{thm}
%
%
The proof of Theorem~\ref{thm:order.nonaut} is presented in Appendix~\ref{app.2}. 
Theorem~\ref{thm:order.nonaut} permits to adjust the already known schemes for autonomous systems for construction of SPRK. Some examples of higher order methods from (\cite{Forest1990}) are
a 3-order symplectic PRK defined by 
\begin{equation} \label{eq:3orderSPRK}
        b = (\frac{7}{24}, \ \frac{3}{4}, \  -\frac{1}{24}), \quad    
        B = (\frac{2}{3},  \ -\frac{2}{3},\ 1);
\end{equation}
and 4-order method given by 
\begin{equation} \label{eq:4orderSPRK}
    \begin{aligned}
        b &= (\frac{7}{48}, \ \frac{3}{8}, \ -\frac{1}{48}, \ -\frac{1}{48}, \ \frac{3}{8}, \ \frac{7}{48}),\\    
        B &= (\frac{1}{3},\ -\frac{1}{3},\ 1,\ -\frac{1}{3},\ \frac{1}{3},\ 0).
    \end{aligned}
\end{equation} 
In numerical tests we also use the following 2-order scheme
\begin{equation} \label{eq:2orderSPRK}
        b = (0, \ 1), \quad    
        B = (\frac{1}{2},  \ \frac{1}{2}).
\end{equation}
\subsection{SPRK network and properties}
Let us now use the constructed higher order SPRK 
for the neural network design. In this case, the activation function is given by an $n$th order explicit SPRK \eqref{eq:expl.SPRK} applied to the Hamiltonian system \eqref{eq:Ham.syst}. This class of neural networks will be called SPRK Net and denoted by
\begin{equation} \label{eq:SPRK.net}    
z_{n+1} = \mathrm{SPRK}(z_n, u_{n,1}, \dots, u_{n,s}),
\end{equation}
for $n = 1, \dots, N-1$ with $z_{n} = (q_{n}, p_{n})$ and $u_{n,i}$ approximating $(u(t_n + c_i h),u(t_n + C_i h)) $ for $i = 1, \dots, s$. 
The symplecticity of the underlying SPRK method directly implies that each step between layers, i.e., from $z_n$ to $z_{n+1}$, satisfies \eqref{eq:sympl.num.preservation}. This implies in particular the non vanishing gradient property which follows from \eqref{eq:sympl.nonvanishing}.
 \begin{cor}
	SPRK Net 
 admits the non vanishing gradient property.
	\end{cor}
The SPRK net has the same structural properties as the Hamiltonain network, which permits to show the universal approximation property of SPRK Net. As in the case of the Hamiltonian net (\cite{Pinkus1999}) the following theorem can be used.
 \begin{thm}[\cite{Pinkus1999}] \label{th:Pinkus.for.UAP}

    Let $\sigma: \R \rightarrow \R$ be a non polynomial continuous function, then for any continuous $f: \Omega \rightarrow \R^n$ with $\Omega \subset \R^n$, any $\varepsilon >0$, there exist $N \in \N, V_j, K_j \in R^n \times R^n$ and $d_j \in \R^n$ such that the function $g : \R^n \rightarrow \R^n$ defined by
    $$g(x) = \sum_{j=0}^{N-1} K_j \sigma( V_j x + d_j)$$
    satisfies
    $$\sup_{x \in \Omega} \| f(x) - g(x) \| \leq \varepsilon.$$
		\end{thm} 
It is possible to represent SPRK net in form of $ g(x) = \sum_{j=0}^{N-1} K_j \sigma( V_j x + d_j)$ for a particular restriction on the parameters in \eqref{eq:network:Hamiltonian}. This leads to the universal approximation property for SPRK Net
\begin{thm} \label{th:univ.approx}
Given a compact $\Omega \in \R^{n}$ and $\xi \in \Omega$, the SPRK network defined by \eqref{eq:SPRK.net} with $(q_0, p_0) = (0, \xi)$, $q_k, p_k \in \R^{n}$ for $k = 0, \dots, N$ and the output $q_N$, that is projection of  $(q_N, p_N)$ to the first component, satisfies the UAP for continuous functions $f : \Omega \rightarrow \Omega$.
		\end{thm} 
 The proof of the theorem follows the lines of  the proof for Hamiltonian Net in (\cite{Zakwan2023}). We give the complete proof in Appendix~\ref{app.1}. The direct application of Theorem~\ref{th:h.Net.UAP2} leads to the general UAP (Definition~\ref{def:UAP}).
 \begin{cor}
    SPRK network defined by \eqref{eq:SPRK.net} with $(q_0, p_0) = (0, \xi)$, $q_k, p_k \in \R^{n}$ for $k = 0, \dots, N$ and the output $q_N$ composed with an appropriate output layer admits the universal approximation property.
 \end{cor}

\subsection{Relations of SPRK to continuous learning setting}
The learning problem based on SPRK Net is defined by \eqref{eq:DL:optimal.problem} with the activation function given by \eqref{eq:SPRK.net}. 
%
%
%
There exists a continuous counterpart of the learning problem in the form of optimal control problem
					\begin{equation} \label{eq:continuous.OCP}
						\begin{aligned}
							& \min_{u} \mathcal{J}(u) = L(\Phi(z(T)), y) + R(u) \\
							&  \dot z(t)= J \nabla H(z(t),u(t)), \\
							&  z(0) = z_0.
						\end{aligned}
					\end{equation}     
					
In this representation \eqref{eq:DL:optimal.problem} is a discretization of \eqref{eq:continuous.OCP}. The order of the approximation will be studied in subsequent publications.

\begin{Rem}
    Notice that there are different ways to discretize the control function $u(t)$ in \eqref{eq:continuous.OCP} to get optimization parameters in \eqref{eq:DL:optimal.problem}. One option is to consider the control parameters exactly at the stages $t_n + c_i h$ and $t_n + C_i h$. In this case, the number of parameters to consider is at least $2s$. Another approach is to use interpolation polynomials to approximate $u(t)$, this approach needs $n$ optimization parameters for $n$th order of approximation, but it requires computations of the interpolation polynomial after each update of parameters in optimization step.
\end{Rem}
%
%
\subsection{Image classification}
We consider the task from (\cite{Galimberti2023}) of points classification. The training and test data are given by blue and red points. For each point, its feature $x_i \in \R^2$ is the position on a 2-dimensional plane and its label $y_i \in \{0,1\}$ is the color of the point. The task is to classify whether a point from the test data is red or blue. 
We compare the Hamiltonian network and SPRK nets of order 2, 3, 4 based on \eqref{eq:3orderSPRK}-\eqref{eq:2orderSPRK} with 12 layers and the same number of optimization parameters for all networks, this is achieved by considering the same parameters at the stages $i=1, \dots, s$ in \eqref{eq:expl.SPRK}. In our implementation we adapt the code by (\cite{Galimberti2023}). We use the features augmentation to dimension 4 and run the optimization with Adams optimizer and 100 epochs. The resulting test accuracy is presented in following table and shows the percentage of the correctly classified points. We can observe very similar results for all the nets with the tendency to be better in the higher order case.  

\begin{tabular}{||c | c| c| c ||} 
 \hline
 H-Net & 2-SPRK Net & 3-SPRK Net &4-SPRK Net \\ [0.5ex] 
 \hline\hline
 63.22 \% & 60.68 \% & 69.65 \% & 73.62 \% \\
 \hline
\end{tabular}
\subsection{Learning Hamiltonian system} \label{sec:example.Hamiltlonan}
The second task is the learning of an autonomous Hamiltonian system. The training data is given by points of trajectories $\{ x(t_i)\}$ generated by numerical integration of an unknown dynamics $\dot{x} = f(t,x)$. The task is to learn $f(t,x)$ in order to generate new trajectories. The considered dynamics is given by the Kepler problem.
\begin{equation} \label{eq:Kepler}
    \dot{q} = p, \quad \dot{p} = - \frac{\pi}{4 \| q\|^{3/2}} \, q
\end{equation}

We assume that the dynamics of $q$ is known and we look for the unknown dynamics $\dot p = f(q)$. The loss function is given by $L = \frac{1}{n}\sum_{i=0}^{n}\| x(t_i) - x_i\|^2$ and we also consider the regularizer of the form $\sum_{n,j} \| u_{n,j}\|^2$. \\
The network construction in this case includes two steps:
\begin{enumerate}
		 	\item approximate $f(q)$ using SPRK-Net such that $z_N \approx f(q)$
		 	\item use the same SPRK method to integrate \eqref{eq:Kepler} 
    to obtain points approximating trajectories  $x_k \approx x(t_k)$ 
		 \end{enumerate}
Based on these two steps, the gradients of the loss with respect to optimization parameters are computed as a composition of gradients from the network approximation of $f(q)$ and the gradients from the numerical integration of \eqref{eq:Kepler} as follows: 
\begin{multline*}
    \sum_{i=0}^{n} \frac{\partial L}{\partial x_i}  \sum_{k = 0}^{i-1} ( \prod_{l = k+1}^{i} \frac{\partial x_{l+1}}{\partial x_{l}})  \frac{\partial x_{k+1}}{\partial  z_{N}}  ( \prod_{s = j+1}^{N-1} \frac{\partial z_{s+1}}{\partial z_{s}}  )\frac{\partial z_{j+1}}{\partial u_j}.  
\end{multline*}
Symplecticity in both parts of the network guarantees the non vanishing gradient property. 

In our implementation, we compare the Hamiltonian net with the 2nd, 3rd and 4th order SPRK nets based on \eqref{eq:3orderSPRK}-\eqref{eq:2orderSPRK}. In all the cases, the same SPRK is used to approximate $f(q)$ and then to integrate the dynamics. The training set is given by 2000 point of the 27 trajectories obtained using sci.integrate.odeint from different initial data and on the same time interval. The test data is given by a trajectory which is generated by the initial condition used for training data, but integrated on a 2-times longer time interval. We use stochastic gradient descent method with 270 epochs. The resulting  test $L_2$-norm accuracy is presented in the following table. The higher order SPRK-Net is more accurate according to the numerical results.

%
\begin{tabular}{||c | c| c| c ||} 
 \hline
 H-Net & 2-SPRK Net & 3-SPRK Net &4-SPRK Net \\ [0.5ex] 
 \hline\hline
 1.26692 & 1.16215 & 0.61231 & 0.18060 \\
 \hline
\end{tabular}

\section{Conclusion and future work}
In this paper, we constructed a class of higher order Hamiltonian nets with the non vanishing gradient and universal approximation properties based on SPRK methods, which generalizes the already existing Hamiltonian Net. The new class of networks has in addition the property to be a higher order approximation of the learning problem in the continuous setting. 
We tested the network on a classification task and learning of a Hamiltonian system. In both learning tasks the higher order network led to a better test accuracy. 

This work opens the doors for different extensions. By construction, SPRK Net allows the learning of separable non autonomous Hamiltonian systems. 
Another possible extension is the construction of higher order reversible networks based on reversible symplectic methods. 
The connection between the new network and the optimal control suggests different parameterization options, which might be advantageous for learning tasks.  


\bibliography{ifacconf}             

\begin{thebibliography}{16}
\providecommand{\natexlab}[1]{#1}
\providecommand{\url}[1]{\texttt{#1}}
\providecommand{\urlprefix}{URL }
\expandafter\ifx\csname urlstyle\endcsname\relax
  \providecommand{\doi}[1]{doi:\discretionary{}{}{}#1}\else
  \providecommand{\doi}{doi:\discretionary{}{}{}\begingroup
  \urlstyle{rm}\Url}\fi

\bibitem[{Abia and Sanz-Serna(1993)}]{abia1993}
Abia, L. and Sanz-Serna, J. (1993).
\newblock Partitioned runge-kutta methods for separable hamiltonian problems.
\newblock \emph{Mathematics of Computation}, 60(202), 617--634.

\bibitem[{Bengio et~al.(1994)Bengio, Simard, and Frasconi}]{Bengio1994}
Bengio, Y., Simard, P., and Frasconi, P. (1994).
\newblock Learning long-term dependencies with gradient descent is difficult.
\newblock \emph{IEEE Transactions on Neural Networks}, 5(2), 157--166.

\bibitem[{Benning et~al.(2021)Benning, Celledoni, Ehrhardt, Owren, and
  Schönlieb}]{Benning2021}
Benning, M., Celledoni, E., Ehrhardt, M.J., Owren, B., and Schönlieb, C.B.
  (2021).
\newblock Deep learning as optimal control problems.
\newblock \emph{IFAC-PapersOnLine}, 54(9), 620--623.
\newblock 24th International Symposium on Mathematical Theory of Networks and
  Systems MTNS 2020.

\bibitem[{Celledoni et~al.(2021)Celledoni, Ehrhardt, Etmann, Mclachlan, Owren,
  Schonlieb, and Sherry}]{CELLEDONI2021}
Celledoni, E., Ehrhardt, M.J., Etmann, C., Mclachlan, R.I., Owren, B.,
  Schonlieb, C.B., and Sherry, F. (2021).
\newblock Structure-preserving deep learning.
\newblock \emph{European Journal of Applied Mathematics}, 32(5), 888–936.

\bibitem[{Chen et~al.(2018)Chen, Rubanova, Bettencourt, and
  Duvenaud}]{Chen2018}
Chen, R.T.Q., Rubanova, Y., Bettencourt, J., and Duvenaud, D.K. (2018).
\newblock Neural ordinary differential equations.
\newblock In S.~Bengio, H.~Wallach, H.~Larochelle, K.~Grauman, N.~Cesa-Bianchi,
  and R.~Garnett (eds.), \emph{Advances in Neural Information Processing
  Systems}, volume~31. Curran Associates, Inc.

\bibitem[{Dupont et~al.(2019)Dupont, Doucet, and Teh}]{Dupont2019}
Dupont, E., Doucet, A., and Teh, Y.W. (2019).
\newblock Augmented neural odes.
\newblock In H.~Wallach, H.~Larochelle, A.~Beygelzimer, F.~dAlch\'{e}-Buc, E.~Fox, and R.~Garnett (eds.), \emph{Advances in Neural
  Information Processing Systems}, volume~32. Curran Associates, Inc.

\bibitem[{Forest and Ruth(1990)}]{Forest1990}
Forest, E. and Ruth, R.D. (1990).
\newblock Fourth-order symplectic integration.
\newblock \emph{Physica D: Nonlinear Phenomena}, 43(1), 105--117.

\bibitem[{Galimberti et~al.(2023)Galimberti, Furieri, Xu, and
  Ferrari-Trecate}]{Galimberti2023}
Galimberti, C.L., Furieri, L., Xu, L., and Ferrari-Trecate, G. (2023).
\newblock Hamiltonian deep neural networks guaranteeing nonvanishing gradients
  by design.
\newblock \emph{IEEE Transactions on Automatic Control}, 68(5), 3155--3162.

\bibitem[{Giesecke and Kröner(2021)}]{Giesecke2021}
Giesecke, E. and Kröner, A. (2021).
\newblock Classification with runge-kutta networks and feature space
  augmentation.
\newblock \emph{Journal of Computational Dynamics}, 8(4), 495--520.

\bibitem[{Hairer et~al.(2013)Hairer, Lubich, and Wanner}]{Hairer2013}
Hairer, E., Lubich, C., and Wanner, G. (2013).
\newblock \emph{Geometric Numerical Integration: Structure-Preserving
  Algorithms for Ordinary Differential Equations}, volume~31.
\newblock Springer Science \& Business Media.

\bibitem[{Kalogiratou et~al.(2011)Kalogiratou, Monovasilis, and
  Simos}]{Kalogiratou2011}
Kalogiratou, Z., Monovasilis, T., and Simos, T.E. (2011).
\newblock Symplectic partitioned runge-kutta methods for the numerical
  integration of periodic and oscillatory problems.
\newblock In T.E. Simos (ed.), \emph{Recent Advances in Computational and
  Applied Mathematics}, 169--208. Springer Netherlands, Dordrecht.

\bibitem[{Matsubara et~al.(2023)Matsubara, Miyatake, and
  Yaguchi}]{Matsubara2023}
Matsubara, T., Miyatake, Y., and Yaguchi, T. (2023).
\newblock The symplectic adjoint method: Memory-efficient backpropagation of
  neural-network-based differential equations.
\newblock \emph{IEEE Transactions on Neural Networks and Learning Systems},
  1--13.

\bibitem[{Pinkus(1999)}]{Pinkus1999}
Pinkus, A. (1999).
\newblock Approximation theory of the mlp model in neural networks.
\newblock \emph{Acta Numerica}, 8, 143–195.

\bibitem[{Sanz-Serna(2016)}]{SanzSerna2016}
Sanz-Serna, J.M. (2016).
\newblock Symplectic runge--kutta schemes for adjoint equations, automatic
  differentiation, optimal control, and more.
\newblock \emph{SIAM Review}, 58(1), 3--33.

\bibitem[{Yan et~al.(2019)Yan, Du, Tan, and Feng}]{yan2019}
Yan, H., Du, J., Tan, V.Y., and Feng, J. (2019).
\newblock On robustness of neural ordinary differential equations.
\newblock \emph{arXiv preprint arXiv:1910.05513}.

\bibitem[{Zakwan et~al.(2023)Zakwan, d’Angelo, and
  Ferrari-Trecate}]{Zakwan2023}
Zakwan, M., d’Angelo, M., and Ferrari-Trecate, G. (2023).
\newblock Universal approximation property of hamiltonian deep neural networks.
\newblock \emph{IEEE Control Systems Letters}, 7, 2689--2694.

\end{thebibliography}
                                                   







\appendix
\section{Proof of Theorem \ref{th:univ.approx}} \label{app.1}

The main idea is to represent the SPRK Net in the form of $ g(x) = \sum_{j=0}^{N-1} K_j \sigma( V_j x + d_j)$ and use Theorem~\ref{th:Pinkus.for.UAP}. This is done using a particular form of parameters in \eqref{eq:H.net}, namely,
	\begin{equation} \label{eq:restricted.params}
	     W(t) = \begin{pmatrix}
		W(t) & 0 \\
		0 & 0
	\end{pmatrix},  \beta(t) = \begin{pmatrix}
		\beta(t) \\
		0
	\end{pmatrix},  \eta(t) = \begin{pmatrix}
		0\\
		-\eta(t).
	\end{pmatrix} 
 \end{equation}
SPRK Net defined by \eqref{eq:SPRK.butcher} then leads to
\begin{align}
  &p_{n+1} = p_n + h \sum_{i=1}^{s} B_i \eta_{i,n}  \label{eq:proof.UAP} \\
  &q_{n+1} = q_n + h \sum_{i=1}^{s} b_i W_{i,n}^\top \sigma(W_{i,n} p_n + h d_{i,n}(\eta) + \beta_{i,n} ), \notag    
\end{align}
where $d_{i,n}(\eta) = \sum_{j=0}^{i-1} B_j W_{i,n} \eta_{j,n}$. 
As the result, the output $q_N$ of the SPRK Net is given by
	$$q_N = \sum_{j=0}^{N-1}  K_j \sigma(  V_j p_0 +   d_j),$$
for a certain matrices $K_j, V_j$ and vector $d_j$ determined from $W_{i,n}, \beta_{i,n}, \eta_{i,n}$ using \eqref{eq:proof.UAP}. Applying Theorem~\ref{th:Pinkus.for.UAP} we obtain the UAP for SPRK Net. If the property holds for the restrictive set of parameters \eqref{eq:restricted.params}, then it also holds for the general class used in the construction of SPRK Net. 

\section{Proof of Theorem \ref{thm:order.nonaut}} \label{app.2}
The main idea of the proof is to show that the new order conditions appearing in case of non autonomous separable system coincide with the order-conditions for autonomous systems if $c_i = \sum_{j = 1}^{s}A_j$ and $C_i = \sum_{j = 1}^{s}a_j$. The proof is based on the comparison of the terms in the Taylor expansion of the analytical and numerical solutions. For simplicity of exposition, we denote $f(p,t) = \frac{\partial}{\partial p} H(q,p,t)$  and $g(q,t) = -\frac{\partial}{\partial q} H(q,p,t)$. 
We write \eqref{general.PRK} as a first order method for $q_{n+1}$
$$q_{n+1} = q_n + h \varphi_q(t_n, q_n, h).$$
At the same time the Taylor expansion for analytic solution $q(t)$ with $q(t_n) = q_n$ and $p(t_n) = p_n$, where all the functions are evaluated at $(q_n, p_n, t_n)$, is given by
    \begin{align}
    &q(t_{n+1}) = q_n + hf + \frac{h^2}{2} \left(f'_p g + f_t\right) +  \label{eq:Taylor.q} \\ 
    &\frac{h^3}{3!} \left(f''_{pp} (g,g) + f'_{p} g'_q f +2f''_{pt}g + f'_pg'_t  + f''_{tt}\right) + \notag \\ 
    &    \frac{h^4}{4!}(f'''_{ppp} (g,g,g) + 3f''_{pp}(g'_qf,g) + f'_pg'_q f'_p g  +  f'_pg''_{qq}(f,f) +
  \notag \\  
    &  3(f'''_{ppt}(g,g) + f''_{pp}(g'_t,g) +f'''_{ptt}g + f''_{pt}g'_q f + f''_{pt}g'_t ) ) + \notag \\
    & 2f'_pg''_{qt}f +   f'_pg'_qf'_t + f'_pg''_{tt} + f''_{ttt})  + O(h^5). \notag
    \end{align}
The same can be written for variable $p$ by exchanging the roles of $f$ and $g$. Therefore, without loss of generality, we can consider the proof for the order of approximation of $q(t)$ only. The consistency condition holds trivially for \eqref{general.PRK}
if $\sum_{i=1}^s b_i = \sum_{i=1}^s B_i = 1$. The second order condition holds if $\left. \frac{\partial}{\partial h}\varphi_q(t_n, q_n, h) \right|_{h=0}$ coincides with second order term in the Taylor expansion of $q$. The method is of the $k$th order if all terms in $\left. \frac{\partial^j}{\partial h^j}\varphi_q(t_n, q_n, h) \right|_{h=0}$ coincide with the $(j+1)$st order terms in the Taylor series for $j<k$.
For the second order condition we consider the first derivative of $\varphi_q$, that is
\begin{multline*}
    \frac{\partial}{\partial h}\varphi_q(t_n, q_n, h) = \sum_{i=1}^s b_i ( f'_p(P_i,t_n + c_i h)\frac{\partial}{\partial h}P_i(h) + \\
    c_i f'_t(P_i,t_n + c_i h)) = \sum_{i,j=1}^s b_iA_{ij}(f_p g + h f'_p(g'_q \frac{\partial}{\partial h}Q_j(h) + C_j g'_t))\\
    + \sum_{i=1}^s b_i c_i f'_t).
\end{multline*}
Therefore, 
\begin{multline*}\left. \frac{\partial}{\partial h}\varphi_q(t_n, q_n, h) \right|_{h=0} = \sum_{i,j=1}^s b_i A_{ij}f_p(p_n,t_n) g(q_n,t_n) \\ + \sum_{i=1}^s b_i c_i f'_t(p_n,t_n)), 
\end{multline*}
where the second term implies new conditions with respect to the condition for non autonomous case. It is easy to see, that $c_i = \sum_{j}A_{ij}$ implies that the new condition $\sum_{i} b_i c_i =  \frac{1}{2}$ coincides with the condition for autonomous systems $\sum_{i,j} b_i A_{ij} =  \frac{1}{2}$. The same can be done for the 3rd order conditions. 
\begin{multline*}
    \frac{\partial^2}{\partial h^2}\varphi_q(t_n, q_n, h)  = \sum_{i,j= 1}^{s} b_iA_{ij}(f''_{pp} (\frac{\partial}{\partial h} P_i(h), \frac{\partial}{\partial h} P_i(h)) \\ + c_i f''_{pt}g  + 2f'_p(g'_q \frac{\partial}{\partial h}Q_j(h) + C_j g'_t))
    + \\ \sum_{i,j= 1}^{s} b_iA_{ij}h \frac{\partial}{\partial h }(f'_p(g'_q \frac{\partial}{\partial h}Q_j(h)  + C_j g'_t)) + \\
    \sum_{i= 1}^{s} b_i c_i (f''_{pt}\frac{\partial}{\partial h}P_i(h) + c_i f''_{tt}).
\end{multline*}
which implies at $h=0$ 
\begin{multline*}
    \left. \frac{\partial^2}{\partial h^2}\varphi_q(t_n, q_n, h) \right|_{h=0} = \sum_{i,j,k= 1}^{s} b_i A_{ij} A_{ik}f''_{pp}(g,g) + \\ 2\sum_{i,j,k= 1}^{s} b_i A_{ij} a_{jk} f'_{p}g'_q f + 2\sum_{i,j= 1}^{s} b_i A_{ij} C_j f'_{p}g'_t  \\ + 2\sum_{i,j= 1}^{s} b_i c_i A_{ij}  f''_{pt}g + \sum_{i= 1}^{s} b_i c_i c_i f''_{tt}.  
\end{multline*}
The first two term are present in the case of the autonomous system and the last three are new and only present in the non autonomous case. The conditions $c_i = \sum_j A_{ij}$ and $C_i = \sum_j a_{ij}$ imply that the order conditions for the autonomous terms and the conditions  for non autonomous terms coincide. Finally let us consider the 4th order conditions.
\begin{multline*}
    \left. \frac{\partial^3}{\partial h^3}\varphi_q(t_n, q_n, h) \right|_{h=0} = \sum_{i,j,k,l} b_i A_{ij} A_{ik} A_{il}f'''_{ppp}(g,g,g) + \\
   3 \sum_{i,k,l} b_i c_i A_{ik} A_{il} f'''_{ppt}(g,g) + 8 \sum_{i,j,k,l} b_i A_{ij} a_{jk} A_{il} f''_{pp}(g'_qf,g) +  \\
   8 \sum_{i,j,l} b_i A_{ij} C_{j} A_{il}f''_{pp}(g'_t,g) + 3 \sum_{i,l} b_i c_{i} c_{i} A_{il}f'''_{ptt}g + \\
   8 \sum_{i,j,k} b_i c_i A_{ij} a_{jk} f''_{pt}g'_qf + 8 \sum_{i,j} b_i c_{i} A_{ij} C_{j} f''_{pt}g'_t + \\
   6 \sum_{i,j,k,l} b_i A_{ij} a_{jk} a_{jl}  f'_{p}g''_{qq}(f,f) + 12 \sum_{i,j,k} b_i A_{ij} a_{jk} C_{j}  f'_{p}g''_{qt}f + \\ 
   12 \sum_{i,j,k,l} b_i A_{ij} a_{jk} A_{kl}  f'_{p}g'_{q}f'_pg+ 12 \sum_{i,j,k} b_i A_{ij} a_{jk} c_{k}  f'_{p}g'_{q}f'_t + \\ 6 \sum_{i,j}b_i A_{ij} C_j C_j f'_{p}g''_{tt} + \sum_i b_i c_i c_i c_i f'''_{ttt}.
    %
\end{multline*}
The terms which appear in the autonomous case are
\begin{multline*}
   \sum_{i,j,k,l} b_i A_{ij} A_{ik} A_{il}f'''_{ppp}(g,g,g) + 8 \sum_{i,j,k,l} b_i A_{ij} a_{jk} A_{il} f''_{pp}(g'_qf,g) +  \\   
   12 \sum_{i,j,k,l} b_i A_{ij} a_{jk} A_{kl}  f'_{p}g'_{q}f'_pg + 6 \sum_{i,j,k,l} b_i A_{ij} a_{jk} a_{jl}  f'_{p}g''_{qq}(f,f). 
    %
\end{multline*}
Comparing the terms with the Taylor expantion of the analytic $q(t)$, we get the conditions
\begin{multline*} b_i A_{ij} A_{ik} A_{il} = \frac{1}{4!}, \ 8\sum_{i,j,k,l} b_i A_{ij} a_{jk} A_{il} =\frac{3}{ 4!}, \\
12 \sum_{i,j,k,l} b_i A_{ij} a_{jk} A_{kl} = \frac{1}{4!}, \ 6 \sum_{i,j,k,l} b_i A_{ij} a_{jk} a_{jl} = \frac{1}{4!}.  
\end{multline*}
And we substitute  $c_i = \sum_j A_{ij}$ and $C_i = \sum_j a_{ij}$ in the new terms appearing in the non autonomous case
\begin{multline*}
   3 \sum_{i,j,k,l} b_i A_{ij} A_{ik} A_{il} f'''_{ppt}(g,g) +   8 \sum_{i,j,k,l} b_i A_{ij} a_{jk} A_{il}f''_{pp}(g'_t,g)  +\\ 3 \sum_{i,j,k,l} b_i A_{ij} A_{ik} A_{il}f'''_{ptt}g +
   8 \sum_{i,j,k,l} b_i A_{ij} a_{jk} A_{il} f''_{pt}g'_qf +\\ 
   8 \sum_{i,j,k,l} b_i A_{ij} a_{jk} A_{il} f''_{pt}g'_t + 
   12 \sum_{i,j,k,l} b_i A_{ij} a_{jk} a_{jl}  f'_{p}g''_{qt}f + \\   12 \sum_{i,j,k,l} b_i A_{ij} a_{jk} A_{kl}  f'_{p}g'_{q}f'_t  +
   6 \sum_{i,j}b_i A_{ij} a_{jk} a_{jl} f'_{p}g''_{tt} + \\
   \sum_{i,j,k,l} b_i A_{ij} A_{ik} A_{il} f'''_{ttt}.
\end{multline*}
The obtained terms coincide with the terms in the Taylor expansion of the analytic solution $q(t)$ under the conditions deduced from the autonomous terms only. This finishes the proof of the theorem for orders $p$ up to 4.
\end{document}